\documentclass[11pt,a4paper]{amsart}
\usepackage{epsfig,a4wide}

\newtheorem{theorem}{Theorem}[section]

\begin{document}

\title{Stick numbers of $2$-bridge knots and links}

\author[Y. Huh]{Youngsik Huh}
\address{Department of Mathematics, School of Natural Sciences, Hanyang University, Seoul 133-791, Korea}
\email{yshuh@hanyang.ac.kr}
\author[S. No]{Sungjong No}
\address{Department of Mathematics, Korea University, 1, Anam-dong, Sungbuk-ku, Seoul 136-701, Korea}
\email{blueface@korea.ac.kr}
\author[S. Oh]{Seungsang Oh}
\address{Department of Mathematics, Korea University, 1, Anam-dong, Sungbuk-ku, Seoul 136-701, Korea}
\email{seungsang@korea.ac.kr}
\thanks{This work was supported by the National Research Foundation of Korea(NRF)
grant funded by the Korea government (MEST) (No.~2009-0074101).}

\begin{abstract}
Negami found an upper bound on the stick number $s(K)$ of a nontrivial knot $K$
in terms of the minimal crossing number $c(K)$ of the knot which is $s(K) \leq 2 c(K)$.
Furthermore McCabe proved $s(K) \leq c(K) + 3$ for a $2$-bridge knot or link,
except in the case of the unlink and the Hopf link.
In this paper we construct any $2$-bridge knot or link $K$ of at least six crossings
by using only $c(K)+2$ straight sticks.
This gives a new upper bound on stick numbers of $2$-bridge knots and links
in terms of crossing numbers.
\end{abstract}

\maketitle

\section{Introduction} \label{sec:intro}

A simple closed curve embedded into the Euclidean $3$-space is called a {\em knot\/}.
A knot $K$ can be embedded many different ways in space.
A {\em stick knot\/} is a knot which consists of finite line segments,
called {\em sticks\/}.
One natural question concerning stick knots may be the {\em stick number\/}
$s(K)$ of a knot $K$ which is defined to be the minimal number of
sticks necessary to construct this stick knot.
Several upper and lower bounds on stick number for various classes of
knots and links were founded.
The most general result is Negami's inequality \cite{N} :
\[ \frac{5+\sqrt{25+8(c(K)-2)}}{2} \leq s(K) \leq 2c(K) \]
for any nontrivial knot or link $K$ other than the Hopf link,
where $c(K)$ is the minimal crossing number of $K$.
Calvo \cite{Ca} improved the lower bound to $\frac{7+\sqrt{8c(K)+1}}{2} \leq s(K)$.
Recently Huh and Oh \cite{HO} utilized the arc index $a(K)$
to determine more precise upper bound,
showing that $s(K) \leq \frac{3}{2} \, (c(K)+1)$ for any nontrivial knot $K$
(especially $s(K) \leq \frac{3}{2} \, c(K)$ for a non-alternating prime knot).
They mainly use the fact that $a(K) \leq c(K)+2$ for any nontrivial knot $K$
in \cite{BP}
and converted each arc presentation of $K$ into a stick knot
by using $\frac{3}{2} \, (c(K)-1)$ sticks.

Precise stick number for specific knots with small crossing number are known
due to work by Randell \cite{R} and Meissen \cite{Me}.
Adams et al. \cite{ABGW} found precise stick numbers for an infinite class of knots,
namely the $(n,n \! \! - \!\! 1)$ torus knots and all their compositions.
In independent work, Jin \cite{J} determined stick numbers for a broader class of knots,
showing that the stick number for any $(p,q)$ torus knot with $p<q<2p$ is $2q$
by using the superbridge index.

There is another way to obtain an upper bound of the stick number
for particular classes of knots and links.
McCabe \cite{Mc} proved that $s(K) \leq c(K) +3$
for any $2$-bridge knot or link $K$ other than the unlink and the Hopf link
by realizing the stick knot visually with $c(K)+3$ sticks.
This upper bound is sharp for small knots and links of at most five crossings.
In \cite{FLS} Furstenberg et al. reduced McCabe's upper bound by $1$
for a few classes of $2$-bridge knots or links, and proposed as an open question
whether or not $s(K) \leq c(K) +2$ for all $2$-bridge knots and links
with crossing number at least $6$.
And Meissen \cite{Me} showed that this inequality holds
for all knots and links with crossing number $7$.
In this paper we give the answer.

\begin{theorem} \label{thm:main}
$s(K) \leq c(K)+2$ for all $2$-bridge knots and links $K$ with $c(K) \geq 6$.
\end{theorem}

\section{Conway notations and integral tangles} \label{sec:conway}

\begin{figure}[h]
\epsfbox{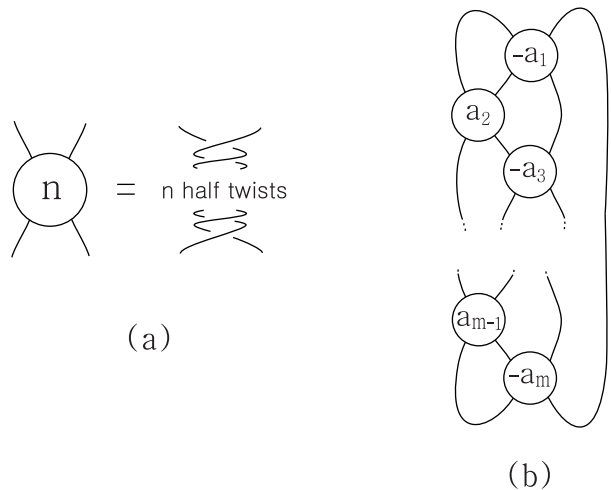}
\caption{Conway notation $(a_1,a_2,\cdots,a_m)$}\label{fig:conway}
\end{figure}

We describe the standard projection of a $2$-bridge knot or link
in terms of the Conway notation which will be useful for the stick construction.
Conway \cite{Co} introduced the concept of tangles, portions of a knot contained
in a topological sphere which intersect the knot exactly four times.
An {\em integral tangle\/} is made from two strands that wrap around each other,
and identified by the number of half-twists (i.e. crossings) within it
as in Figure \ref{fig:conway}(a).
More precisely the integer inside the circle is positive
if it indicates the number of right-handed half-twists and negative if left-handed.
These integral tangles are connected together as in Figure \ref{fig:conway}(b)
to form a $2$-bridge knot or link which is represented by a {\em Conway notation\/}
$(a_1,a_2,\cdots,a_m)$.
Note that if all $a_i$ are positive integers,
then the positive and negative signs of integers in the figure appear alternately,
so it gives a non-nugatory alternating projection of a $2$-bridge knot or link.

As summarized in \cite[Section 2]{Mc}, these are well known facts that
any $2$-bridge knot or link can be represented by Conway notation
$(a_1,a_2,\cdots,a_m)$ with positive integers $a_i$ and odd number $m$
due to work by Burde and Zieschang \cite{BZ},
and this non-nugatory alternating projection displays the minimal number of crossings
due to Kauffman \cite{K}, Murasugi \cite{Mu}, and Thistlethwaite \cite{T}.

To simplify the cases of the main proof we may assume that $a_m \geq 2$.
For, when $a_1 \geq 2$ and $a_m = 1$, we use the Conway notation $(a_m,a_{m-1},\cdots,a_1)$
instead of the original one, both of which indicate the same knot or link.
And when $a_1 = a_m = 1$, we can use its mirror image $(a_2+1,a_3,\cdots,a_{m-2},a_{m-1}+1)$
which has eventually the same crossing number and stick number as the original one.

\begin{figure}[h]
\epsfbox{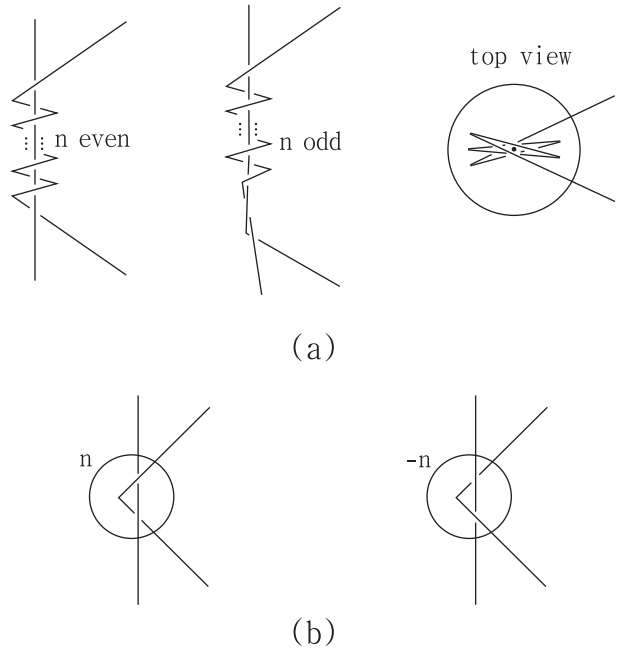}
\caption{$n$-tangles}\label{fig:tangle}
\end{figure}

In the rest of this section,
we illustrate how to construct integral $\pm n$-tangles with $n \geq 2$ by using $n+1$ sticks.
Now observe the $n$-tangle which consists of right-handed $n$ half-twists
as in Figure \ref{fig:tangle}(a).
In this construction, this tangle with $n$ even has one vertical stick,
called a {\em core stick\/}, and the remaining $n$ sticks wrapping around the core stick.
And the tangle with $n$ odd has two sticks which are almost parallel and very close
to each other at their one endpoints so that they look like two slightly perturbed
broken pieces of the core stick,
and $n-1$ sticks wrapping around these two sticks as in the figure.
Both tangles have right-handed $n$ half-twists respectively, and so $n$ crossings.
The right figure shows the drawing of the tangle with $n$ even in the view point from the top.
The $-n$-tangle is simply the mirror image of the positive one.
We use symbols with circles and numbers as drawn in Figure \ref{fig:tangle}(b)
to indicate these tangles.

\section{Proof of Theorem \ref{thm:main}} \label{sec:main}

Let $K$ be any $2$-bridge knot or link with $c(K) \geq 6$.
Recall that $K$ has a Conway notation $(a_1,a_2,\cdots,a_m)$
where $m$ is odd and all $a_i$ are positive integers with $a_m \geq 2$.
Here $c(K) = a_1 + a_2 + \cdots + a_m$.

\begin{figure}[h]
\epsfbox{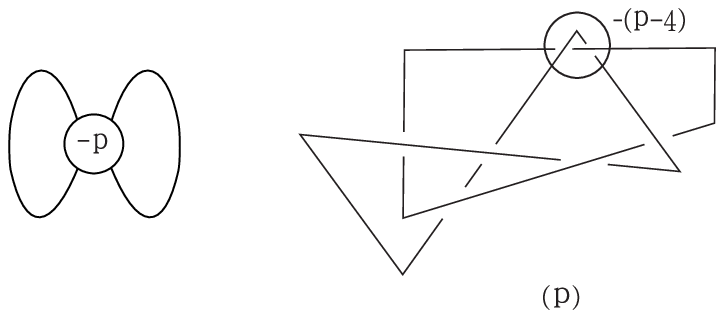}
\caption{$(p)$ case}\label{fig:p}
\end{figure}

\begin{figure}[h]
\epsfig{file=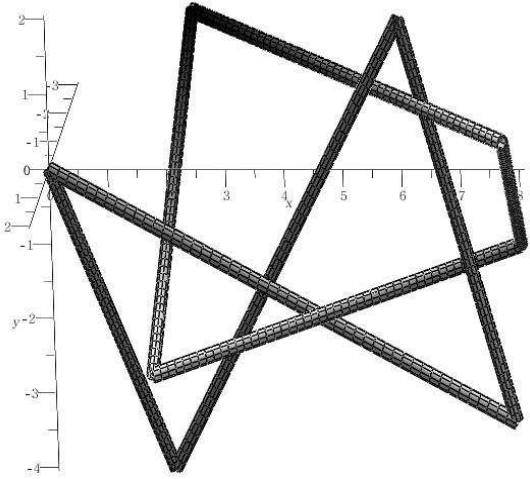, width=8cm}
\caption{visualization of $(6)$ case}\label{fig:6}
\end{figure}

First we construct $K$ of a Conway notation $(p)$ when $p \geq 6$
by using $p+2$ sticks.
$K$ can be obtained from the $-p$-tangle by closing off the four ends
in the standard way.
As drawn in Figure \ref{fig:p},
we utilize the $-(p-4)$-tangle consisting of $p-3$ sticks
and add the remaining $5$ sticks to make $4$ more left-handed half-twists as needed.
Figure \ref{fig:6} shows a specific realization of the link $(6)$
which is the union of two unknotted circles both of which consist of four sticks.
The coordinates of vertices are $\{(0,0,0),(2,-4,0),(6,2,0),(8,-3,1)\}$
and $\{(2,-2,2),(8,-1,0),(8,1,1.7),(2,1,-3)\}$.
Following the tangle construction in the previous section
we can replace the $-2$-tangle by the $-(p-4)$-tangle.
Hence this specific example verifies that our scheme is also realizable for any case $p>6$.

\begin{figure}[h]
\epsfbox{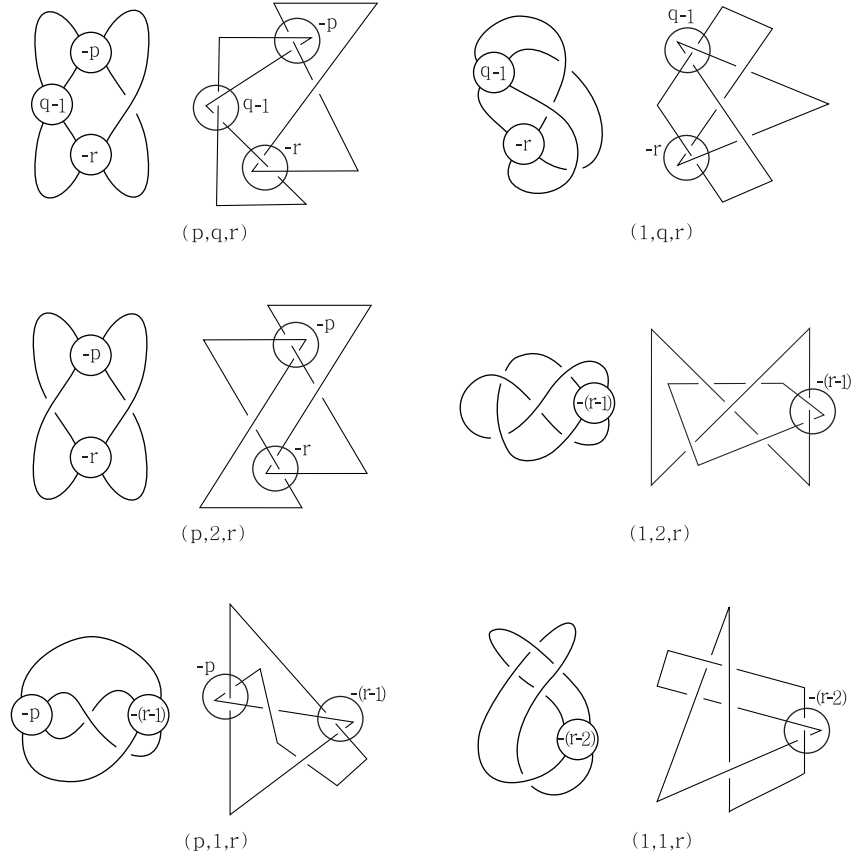}
\caption{$(p,q,r)$ case}\label{fig:pqr}
\end{figure}

Now we construct $K$ of a Conway notation $(p,q,r)$
for positive integers $p,q,r$ with $r \geq 2$ by using $p+q+r+2$ sticks.
Assume that $p \geq 2$ and $q \geq 3$ for the general case.
As in the first link in Figure \ref{fig:pqr},
we rotate the $-p$-tangle part around a horizontal line
for untwisting one half-twist of the $q$-tangle.
Since all three integers $p, q-1$, and $r$ are at least $2$,
we can construct these $-p$-tangle, $(q-1)$-tangle and $-r$-tangle
by using $p+q+r+2$ sticks in total
as described in the last paragraph of Section \ref{sec:conway}.
Now connect them together to build a knot or link $(p,q,r)$ as drawn in the figure
so that two pairs of sticks are joined into two sticks
and two new sticks are added for closing off.
Thus the total number of sticks are unchanged as desired.

Easy calculation tells us that there are five more non-general cases as follows;
$(p,2,r)$ with $p \geq 2$ and $r \geq 2$, $(p,1,r)$ with $p \geq 2$ and $r \geq 3$,
$(1,q,r)$ with $q \geq 3$ and $r \geq 2$, $(1,2,r)$ with $r \geq 3$,
and $(1,1,r)$ with $r \geq 4$.
Note that the case $(p,1,r)$ with $p \geq 3$ and $r \geq 2$ has the same result
as $(p,1,r)$ with $p \geq 2$ and $r \geq 3$.
Similar constructions to the general case can be applied for these five cases.
The results of all the cases are drawn in Figure \ref{fig:pqr}.
One can see that each has exactly $p+q+r+2$ sticks.
To guarantee that all these six cases can be realized in $3$-space as intended,
we provide the vertex-coordinates of the simplest case of each,
such as $(2,3,2)$, $(1,3,2)$, $(2,2,2)$, $(1,2,3)$, $(2,1,3)$ and $(1,1,4)$ in the following table.
Note that each of the cases $(2,3,2)$, $(2,2,2)$ and $(1,2,3)$ consists of two circles,
so they need two sets of vertices indicating two circles respectively. \\

{\small
\begin{tabular}{|c|c|}    \hline
$(p,q,r)$   &   coordinates of vertices   \\ \hline
$(2,3,2)$   &   \{(0,0,0),(0,10,0),(5,10,0),(-2,5,1),(5,0,-4)\}, \{(3,11,-1),(7,2,5),(1,2,-2),(7,11,5)\}  \\
$(1,3,2)$   &   \{(0,0,7),(3,-5,-8),(4,-4,2),(1,3,-1.3),(5,0,0),(1,-3,0),(4,4,0),(3,5,-10)\}  \\
$(2,2,2)$   &   \{(0,0,0),(4,0,0),(0,-5,0),(4,-5,-1)\}, \{(2,-4,-2),(6,-4,3),(2,1,-1),(6,1,5)\}  \\
$(1,2,3)$   &   \{(0,0,0),(6,-5,0),(6,0,0),(0,-5,1)\}, \{(1.5,-5,3),(7,-2.5,-1),(5.5,-1.5,1),(1,-1.5,-2.5)\}  \\
$(2,1,3)$   &   \{(0,0,0),(0,-10,0),(9,-5,-1),(-2,-5,6),(4,-4,-15),(4,-6.5,35),(9,-8,-130),(10,-7,0)\}  \\
$(1,1,4)$   &   \{(1,2.5,2.1),(7,2,-5),(7,-5,5),(3,-5,0),(3,4,0),(1,-4,0),(8,-2,1),(0,0,-1)\} \\  \hline
\end{tabular}
}
\\ \\

\begin{figure}[h]
\epsfbox{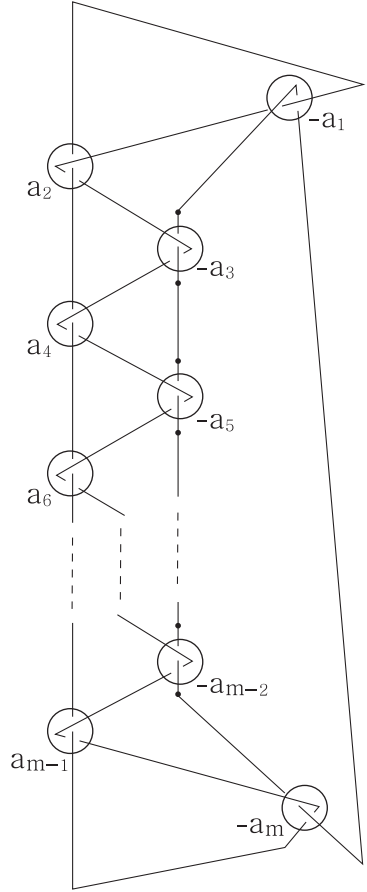}
\caption{general case}\label{fig:general}
\end{figure}

Finally we construct $K$ of a Conway notation $(a_1,a_2,\cdots,a_m)$
for the cases of $m \geq 5$ by using $a_1 + a_2 + \cdots + a_m+2$ sticks,
that is $c(K) + 2$ sticks.
First we consider the general cases that all $a_i$ are greater than $1$.
See Figure \ref{fig:general}.
We construct all $\pm a_i$-tangles individually
by using $a_1 + a_2 + \cdots + a_m+m$ sticks in total.
Now connect them together to build the knot or link as drawn in the figure.
In this scheme, all $a_i$-tangles with $i$ even share just one core stick,
and each pair of $\pm a_i$-tangle and $\pm a_{i+1}$-tangle for $i=1, \cdots,m-1$
share exactly one stick.
On the other hand we have to add some extra sticks to connect
each pair of $\pm a_{2i+1}$-tangle and $\pm a_{2i+3}$-tangle
for $i=1, \cdots,\frac{m-5}{2}$,
and two more sticks to connect the left-most core stick to
$-a_1$-tangle and $-a_m$-tangle.
Thus the total number of sticks is
\[ a_1 + a_2 + \cdots + a_m+m-((\frac{m-1}{2} - 1) + (m-1)) + (\frac{m-5}{2} + 2) \]
that is $a_1 + a_2 + \cdots + a_m+2$.

\begin{figure}[h]
\epsfbox{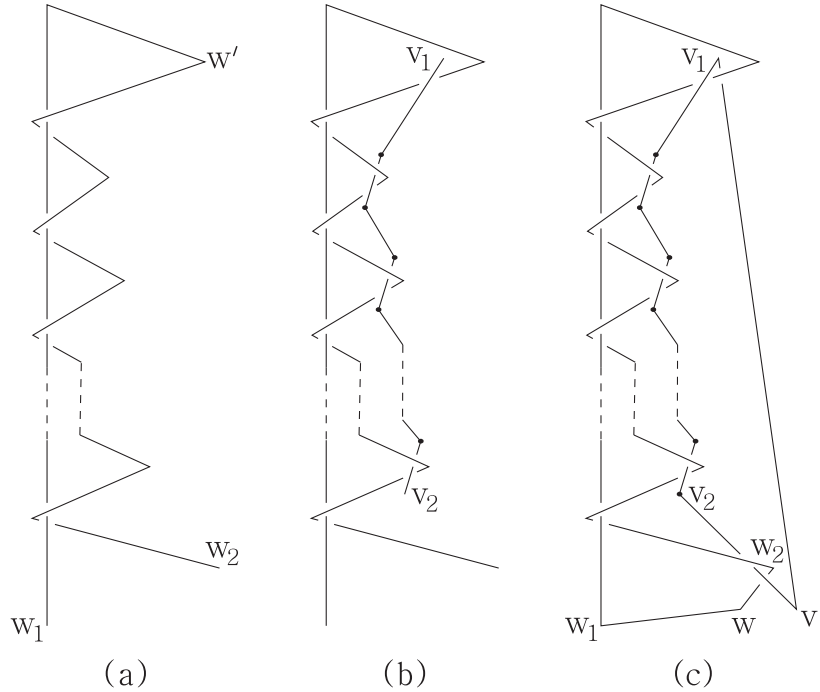}
\caption{realization of the general case}\label{fig:realization}
\end{figure}

Now we illustrate how to build up a stick presentation of a special case
where all $a_i$'s are $2$ in three stages as in Figure \ref{fig:realization}.
Note that one can easily apply this construction to general cases
where some of $a_i$'s are greater than $2$.
At the first stage of this construction,
Figure \ref{fig:realization}(a) shows how to build the main frame of this stick presentation.
Starting from the origin point $w_1=(0, 0, 0)$ in $\mathbb{R}^3$,
take a long vertical stick along $y$-axis,
and more sticks wrapping around this vertical stick from top to bottom
so that these non-vertical sticks rotate clockwise in the view point from the top.
Let $w_2$ be the terminal point of this part.
The intermediate point $w^{\prime}$ in the figure is placed on the $xy$-plane
so that its $x$-coordinate is positive.
In Figure \ref{fig:realization}(b) we see how a connected series of sticks runs back and forth
between the sticks constructed already.
The two end points of this connected series of sticks are denoted by $v_1$ and $v_2$.
Finally as in Figure \ref{fig:realization}(c)
we connect $v_1$ and $v_2$ by two sticks which share one end point $v$.
Obviously the vertex $v$ must lie far below the $xy$-plane.
Also we connect $w_1$ and $w_2$ by two sticks sharing one end point $w$
which lies farther below than the the vertex $v$ with respect to the $z$-coordinate.

\begin{figure}[h]
\epsfbox{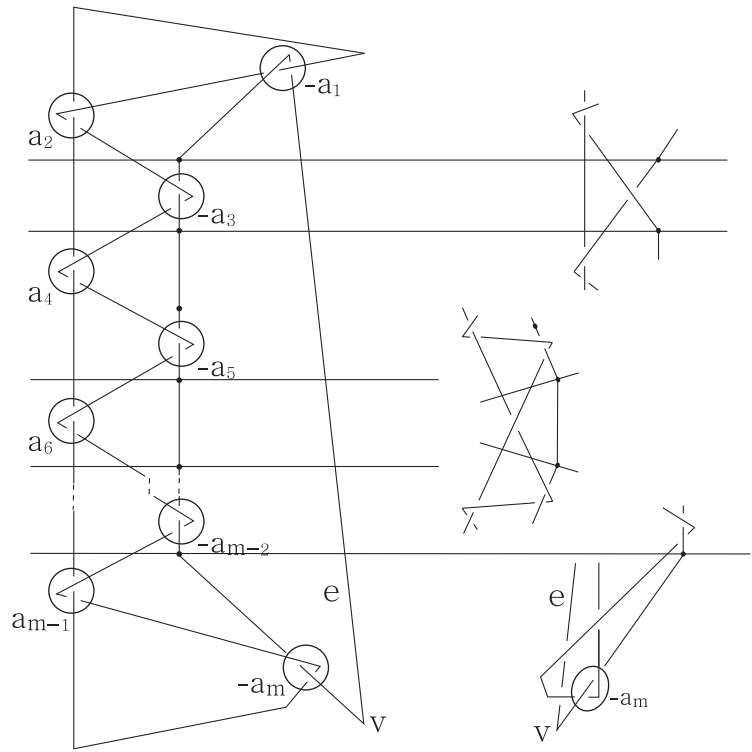}
\caption{some $a_i$'s are $1$}\label{fig:1}
\end{figure}

\begin{figure}[h]
\epsfbox{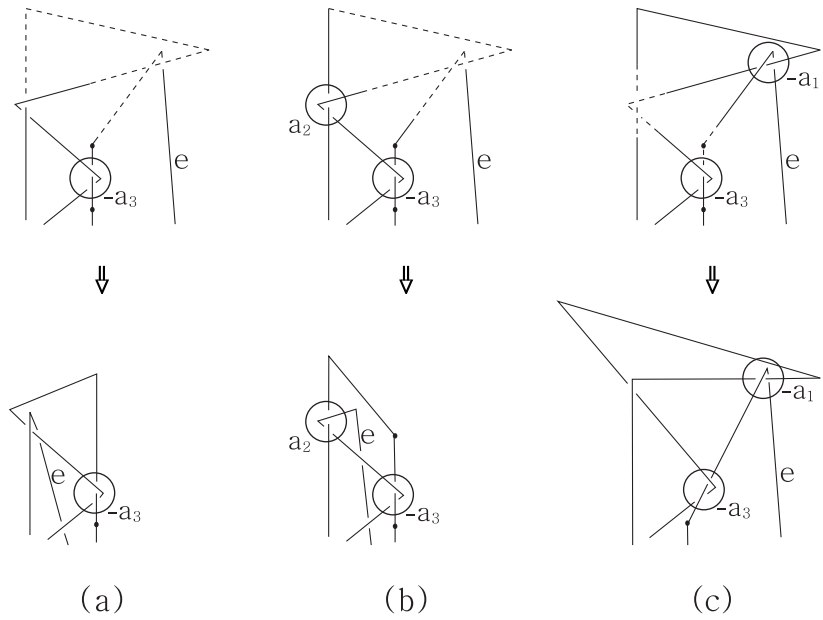}
\caption{$a_1$ or $a_2$ are $1$}\label{fig:11}
\end{figure}

Now we handle the cases when $a_i$'s are $1$.
We distinguish six cases.
The first three cases are illustrated in Figure \ref{fig:1}
and the other cases are in Figure \ref{fig:11}.
First consider the case $a_i=1$ for some $i$ among odd integers $3,5,\cdots,m-2$.
Replace the subgraph of the main graph between two horizontal lines
by the graph indicating $a_i=1$ drawn on the right top in Figure \ref{fig:1}.
Next consider the case $a_i=1$ for some $i$ among even integers $4,6,\cdots,m-3$.
Replace the related subgraph by the graph drawn on the right middle.
In this case the main frame will be bent slightly.
For the case $a_{m-1}=1$, replace the related subgraph by the graph on the right bottom.
Note that the stick $e$ does not pass the other sticks
because the vertex $v$ lies far below the $xy$-plane.
Now consider the next three cases where both or either $a_1$ or $a_2$ is $1$.
Figure \ref{fig:11}(a) illustrates how to realize our scheme for the case $a_1=a_2=1$.
At the top of the main graph we remove the dotted line segments
and connect the two pairs of end points properly.
Note again that the stick $e$ does not pass the other sticks.
The last two cases when only $a_1=1$ and only $a_2=1$ can be constructed following the procedure
in Figure \ref{fig:11}(b) and (c) respectively.

By counting the number of sticks in each case,
we see that this stick presentation yields the inequality of the main theorem.

\end{document}